\documentclass[11pt,en]{elegantpaper}
\usepackage{latexsym,color,amsmath,amsthm,amssymb,amscd,amsfonts}
\usepackage{lipsum}

\theoremstyle{plain}
\newtheorem{thm}{Theorem}[section]
\newtheorem{lem}[thm]{Lemma}

\newtheorem{conj}[thm]{Conjecture}

\def\r{\mathbb{R}^n}
\def\s{\mathbb{S}^{n-1}}
\def\k{\mathcal{K}^n_o}

\def\h{H_{u,K}^-}

\title{Two new proofs of partial Godbersen's Conjecture}
\author{Lin Cheng}
\usepackage{array}

\begin{document}

\maketitle

\begin{abstract}
Two new proofs are provided, offering two new perspectives on Godbersen's conjecture when $j=1$ or $n-1$. One of the proofs utilizes Helly's theorem to provide a concise and elegant proof of the inequality in Godbersen's conjecture. The other proof utilizes the Brunn-Minkowski inequality to provide a completely new proof of the inclusion $-K\subset nK$ for convex bodies $K$ with centroid at the origin, thereby proving Godbersen's conjecture.
\end{abstract}

\section{Introduction}

In this article we investigate the new proofs of Godbersen’s conjecture, which was suggested in 1938 by Godbersen \cite{godbersen1938satz} (and independently by Makai Jr. \cite{hajnal1976research}).

\begin{conj}[Godbersen's conjecture]\label{conj:god}
For any convex body $K\subset \mathbb{R}^n$ and any $1\le j\le n-1$,  
\begin{equation}\label{eq:Godbersen-conj} V(K[j], -K[n-j])\le \binom{n}{j} V(K),\end{equation}
with equality holds if and only if $K$ is a simplex. 
\end{conj}

The cases $j=1$ and $j=n-1$ of Conjecture \ref{conj:god} follow from the fact that $-K\subset nK$ for convex body $K$ whose centroid is at the origin (see \cite{bonnesen1934theorie}, page 53), and inclusion which is tight for the simplex \cite{schneider2009stability}.

\begin{thm}\label{1}
For any convex body $K\subset \mathbb{R}^n$ and $j=1$ or $j=n-1$,  
\begin{equation*} V(K[j], -K[n-j])\le n V(K),\end{equation*}
with equality holds if and only if $K$ is a simplex. 
\end{thm}

The other cases are only verified for special convex bodies, such as simplices (which are the equality case) and convex bodies of constant width, as shown in \cite{godbersen1938satz}. Moreover, this fact gives the bound 
$$ V(K[j], -K[n-j])\le n^{\min\{j,n-j\}}V(K),\enspace \text{for}\enspace 1 \leq j \leq n-1.$$
Recently, the paper \cite{artstein2015godbersen}
shows that for any $\lambda \in [0,1]$ and for any convex body $K$ one has that 
$$
\lambda^j (1-\lambda)^{n-j} V(K[j], -K[n-j])\le {V(K)}.
$$
In particular, picking $\lambda = \frac{j}{n}$, we get that 
$$
V(K[j], -K[n-k])\le \frac{n^n}{j^j (n-j)^{n-j}}V(K)\sim \binom{n}{j} \sqrt{2\pi \frac{j(n-j)}{n}}. 
$$

Back to Theorem \ref{1}, this article is organized as follows. In Section 2, some basic facts on convex geometry are showed. In Section 3, a combinatorial approach to Theorem \ref{1} is introduced. Helly's theorem is used to reduce the general case to the case when $K$ is a simplex. In Section 4, Theorem \ref{1} is proved by a geometric inequality for a specific class of concave functions, and the Brunn-Minkowski inequality is used to connect convex bodies and concave functions.

\section{Preliminaries}
The setting for this article is the n-dimensional Euclidean space, $\mathbb{R}^n$. A convex body
is a compact convex set that has a nonempty interior. Denote by $\mathcal{K}^n_o$ the set of convex bodies in $\mathbb{R}^n$ with the origin $o$ in their interiors. A polytope in $\mathbb{R}^n$ is the convex hull
of a finite set of points in $\mathbb{R}^n$ provided it has positive volume $V_n$ (i.e., $n$-dimensional
Lebesgue measure). If the dimension is clear, we write $V_n$ as $V$. Write $\mathcal{P}_o^n$ for the set of polytopes in $\mathbb{R}^n$ with the origin
in their interiors.

The standard inner product of the vectors $x, y \in \mathbb{R}^n$ is denoted by $x\cdot y$. We write
 $\mathbb{S}^{n-1} = \{x \in \mathbb{R}^n:| x |= 1\}$ for  unit sphere in $\mathbb{R}^n$. The letter $\mu$ will be used exclusively to denote a finite Borel measure on $\mathbb{S}^{n-1}$.
For such a measure $\mu$, we denote by supp$\mu$ its support set.

The support function $h_K : \r \to \mathbb{R}$ of a convex body $K$ is defined, for $x \in \r$, by
$$
h_K(x) = \max\{x \cdot y : y \in K\}.
$$
Observe that support functions are positively homogeneous of degree one and subadditive.
The set $\k$ is often equipped with the Hausdorff metric $\delta$. For $K,L \in \k$,
$$
\delta(K, L) = \sup_{u\in \s}
|h_K(u)-h_L(u)|.
$$
In particular, $\mathcal{P}_o^n$ is a dense subset of $\k$ with the Hausdorff metric.

A hyperplane of $\r$ can be written in the form
$$
H_{u,\alpha} = \{x \in \r : x \cdot u = \alpha\}
$$
with $u \in \r \backslash \{o\}$ and $\alpha \in \mathbb{R}$. The hyperplane $H^-_{u,\alpha}$ bounds a closed halfspace
$$
H_{u,\alpha}^- = \{x \in \r : x \cdot u \leq \alpha\}.
$$

Recall that for convex bodies $K_1, \ldots,K_m \subset {\mathbb R}^n$, and non-negative real numbers $\lambda_1, \ldots,\lambda_m$, the volume of $\lambda_1 K_1 + \cdots + \lambda_m K_m$ is a homogeneous nth degree polynomial in the
$\lambda_1,\ldots, \lambda_m$,
\begin{equation*}\label{Eq_Vol-Def}
V\left(\sum_{i=1}^m \lambda_i K_i\right) = 
\sum_{i_1,\dots,i_n=1}^m \lambda_{i_1}\cdots\lambda_{i_n} V(K_{i_1},\dots,K_{i_n}),
\end{equation*}
and the coefficients $V(K_{i_1},\dots,K_{i_n})$,  called the mixed volume of $K_{i_1}, \ldots, K_{i_n}$, are nonnegative, symmetric in the indices, translation invariant and  dependent only on $K_{i_1}, \ldots, K_{i_n}$. $V(K[j], T[n-j])$ denotes the mixed volume of $j$ copies of the convex body $K$ and $n-j$ copies of the convex body $T$.

The surface area measure $S_K$ of a convex body $K$ is a finite Borel measure on $\s$, defined for every Borel set $\omega\subset \s$ by
$$
S_K(\omega)=\mathcal{H}^{n-1}(\nu_K^{-1}(\omega)),
$$
where $\nu_K:\partial K \to \s$ is the Gauss map of $K$ and $\mathcal{H}^{n-1}$ denotes the $(n-1)$-dimensional Hausdorff measure. Moreover, for convex bodies $K$ and $T$, 
\begin{equation}\label{2}
    V(K[1],T[n-1])=\frac{1}{n}\int_{\s}h_K(u)dS_T(u).
\end{equation}
More details could be found in \cite{schneider2014convex}.
\section{From simplex to the general case}
Because of equation (\ref{2}) and that mixed volume is translation invariant, a natural way to consider Theorem 1.2 is to ask whether there is a point $a\in \r$ such that
\begin{equation}\label{3}
h_{-K+a}(u)\leq n h_{K-a}(u)
\end{equation}
for any $u\in$ supp$S_K$. Moreover, equation (\ref{3}) is equivalent to
\begin{equation}\label{4}
    a\cdot u \leq \frac{n}{n+1}h_K(u)-\frac{1}{n+1}h_K(-u).
\end{equation}
For convenience, $H_{u,\frac{n}{n+1}h_K(u)-\frac{1}{n+1}h_K(-u)}^-$ is denoted by $H_{u,K}^-$ and denote $\cap_{u\in \text{supp}S_K} \h$ by $A_K$. 

If $A_K\ne \emptyset$, for $a\in A_K$, equation (\ref{3}) is right for $u \in \text{supp}S_K$ and 
\begin{align*}
    V(-K[1],K[n-1])&=V(-K+a[1],K[n-1])
    \\
    &=\frac{1}{n}\int_{\s}h_{-K+a}(u)dS_K(u)
    \\
    &\leq\int_{\s}h_{K-a}(u)dS_K(u)
    \\
    &=nV(K).
\end{align*}

Therefore, we are going to prove the following theorem in fact.
\begin{thm}\label{3.1}
    For any convex body $K\subset\r$, $A_K\ne \emptyset$.
\end{thm}
Before proving Theorem \ref{3.1}, some essential lemmas are required.
\begin{lem}\label{3.2}
    For any convex body $K\subset\r$ and any $\phi\in GL_n(\r)$, $A_K\ne \emptyset$ is equivalent to $A_{\phi K}\ne \emptyset$.
\end{lem}
\begin{proof}
    According to the definition of support function and surface area measure, 
    \begin{align*}
        A_K\ne \emptyset. &\Longleftrightarrow \cap_{u\in \text{supp}S_K} \h \ne \emptyset.
        \\
        &\Longleftrightarrow \phi(\cap_{u\in \text{supp}S_K} \h)\ne \emptyset.
        \\
        &\Longleftrightarrow \cap_{u\in \text{supp}S_{\phi K}} H_{u,\phi K}^- \ne \emptyset.
        \\
        &\Longleftrightarrow A_{\phi K}\ne \emptyset.
    \end{align*}
\end{proof}

\begin{lem}\label{3.3}
    If $K$ is a simplex in $\r$, then $A_K$ is a one point set.
\end{lem}
\begin{proof}
    According to Lemma \ref{3.2}, it suffices to show that $A_K$ is a one point set if $K$'s vertices are precisely
    the origin $o$ and points $(1,0,\ldots,0),(0,1,\ldots,0),\ldots,(0,0,\ldots,1)$. By direct calculation, 
    $$A_K=\left\{\left(\frac{1}{n+1},\frac{1}{n+1},\ldots,\frac{1}{n+1}\right)\right\},$$
    which means that $A_K=\{\text{centroid of}\thinspace K\}$ if $K$ is a simplex.
\end{proof}
The next theorem is the key to Theorem \ref{3.1}.
\begin{thm}[Helly's theorem \cite{helly1923mengen}]Let $\mathcal{A}$ be a family of at least $n+1$ compact convex sets in $\r$ and assume that any $n+1$ sets in $\mathcal{A}$ have a nonempty intersection. Then, there is a point $x\in \r$ which is contained in all sets of $\mathcal{A}$.
\end{thm}

After all these preparations, now we can prove Theorem \ref{3.1}.
\begin{proof}[\rm{\textbf{Proof of Theorem} \ref{3.1}}]
    According to Helly's theorem, it suffices to show that 
    $$
    \cap_{i=1}^{n+1} H_{u_i,K}^-\ne \emptyset
    $$
    for any different $u_1,\ldots,u_{n+1} \in \text{supp}S_K$. Without loss of generality, assume that $K\in \k$. We prove this theorem by induction on $n$.

    When $n=2$, according to Helly's theorem, it suffices to show that 
    $$
    \cap_{i=1}^{3} H_{u_i,K}^-\ne \emptyset
    $$
    for any different $u_1,u_2,u_3 \in \text{supp}S_K$. Since the rank of $\{u_1,u_2,u_3\}$ is $2$, there exists $\phi \in GL_2(\mathbb{R}^2)$ such that $\{\phi(u_i),\phi(u_j)\}$ form an orthogonal basis of $\mathbb{R}^2$ for some $1\leq i<j\leq 3$. Similarly to the proof of Lemma \ref{3.2}, 
    $$
    \cap_{i=1}^{3} H_{u_i,K}^-\ne \emptyset \Longleftrightarrow
    \cap_{i=1}^{3} H_{\phi(u_i),\phi^{-T}(K)}^-\ne \emptyset.
    $$
    Without loss of generality, assume that $\{u_1,u_2\}$ is an orthogonal basis of $\mathbb{R}^2$. Thus there exist $b_1,b_2\in \mathbb{R}$ such that $$
    u_3=b_1u_1+b_2u_2.
    $$
    Without loss of generality, let $b_1\leq b_2$. Moreover $\cap_{i=1}^{3} H_{u_i,K}^-\ne \emptyset$ is equivalent to that there exist $a_1,a_2\in \mathbb{R}$ such that
    \begin{equation}\label{5}
        \begin{aligned}
        a_1&\leq \frac{2}{3}h_K(u_1)-\frac{1}{3}h_K(-u_1),
        \\
        a_2&\leq \frac{2}{3}h_K(u_2)-\frac{1}{3}h_K(-u_2),
        \\
        a_1b_1+a_2b_2&\leq \frac{2}{3}h_K(u_3)-\frac{1}{3}h_K(-u_3).
    \end{aligned}
    \end{equation}
    
    If $b_2>0$, there always exist $a_1$ and $N \in \mathbb{Z}$ such that for every $a_2\geq N$ the inequality (\ref{5}) holds. 
    
    If $b_2=0$, then $b_1=-1$ since $u_3\in \text{supp}S_K$. Thus inequality (\ref{5}) turns into
    \begin{equation}\label{6}
        \begin{aligned}
        \frac{1}{3}h_K(u_1)-\frac{2}{3}h_K(-u_1)\leq a_1&\leq \frac{2}{3}h_K(u_1)-\frac{1}{3}h_K(-u_1),
        \\
        a_2&\leq \frac{2}{3}h_K(u_2)-\frac{1}{3}h_K(-u_2).
    \end{aligned}
    \end{equation}
    Notice that $o\in K$ and $h_K(u)\geq 0$ for $u\in \mathbb{S}^1$, such $a_1,a_2$ always exist. 
    
    If $b_2<0$, denote $\cap_{i=1}^{3} H_{u_i,h_K(u_i)}^-$ by $L_2$. In particular, $L_2$ is a simplex with $K\subset L_2$ and $A_{L_2}\ne \emptyset$ according to Lemma \ref{3.3}. Moreover, 
    \begin{align*}
        h_K(-u_i)\leq h_{L_2}(-u_i)\enspace \text{and}\enspace h_K(u_i)=h_{L_2}(u_i)
    \end{align*}
    for $i=1,2,3$. Thus $A_{L_2}\subset \cap_{i=1}^{3} H_{u_i,K}^-$ and $\cap_{i=1}^{3} H_{u_i,K}^-\ne \emptyset$. Therefore $A_K\ne \emptyset$ and Theorem \ref{3.1} is right when $n=2$.

    Assume that the case when $n=k-1$ is right. When $n=k$, according to Helly's theorem, it suffices to show that 
    $$
    \cap_{i=1}^{k+1} H_{u_i,K}^-\ne \emptyset
    $$
    for any different $u_1,\ldots,u_{k+1} \in \text{supp}S_K$.
    If $\text{rank}\{u_1,\ldots,u_{k+1}\}<k$, there exists $u_0\in \mathbb{S}^k$ such that $u_0\cdot u_i=0$ for every $i=1,\ldots,k+1$. Consider $P_{u_0^{\perp}}(K)$ as a $(k-1)$-dimensional convex body and notice that 
    $$
    h_K(u_i)=h_{P_{u_0^{\perp}}(K)}(u_i)\enspace \text{and} \enspace h_K(-u_i)=h_{P_{u_0^{\perp}}(K)}(-u_i)
    $$
    for $i=1,\ldots,k+1$. Thus we have $A_{P_{u_0^{\perp}}(K)}\ne \emptyset$ by induction and $\cap_{i=1}^{k+1} H_{u_i,K}^-\ne \emptyset$ since 
    $$
    \frac{k}{k+1}>\frac{k-1}{k}\enspace \text{and}\enspace \frac{1}{k+1}<\frac{1}{k}.
    $$
    If $\text{rank}\{u_1,\ldots,u_{k+1}\}=k$, without loss of generality, assume that $\{u_1,\ldots,u_k\}$ is an orthogonal basis of $\mathbb{R}^k$, and
    $$
    u_{k+1}=b_1u_1+\cdots+b_k u_k
    $$
    with $b_1\leq\cdots\leq b_{k}$. $\cap_{i=1}^{k+1} H_{u_i,K}^-\ne \emptyset$ is equivalent to that there exist $a_1,\ldots,a_k\in \mathbb{R}$ such that
    \begin{equation}\label{7}
        \begin{aligned}
        a_1&\leq \frac{k}{k+1}h_K(u_1)-\frac{1}{k+1}h_K(-u_1),
        \\
        a_2&\leq \frac{k}{k+1}h_K(u_2)-\frac{1}{k+1}h_K(-u_2),
        \\
        &\vdots
        \\
        a_k&\leq \frac{k}{k+1}h_K(u_k)-\frac{1}{k+1}h_K(-u_k),
        \\
        a_1b_1+\cdots+a_kb_k&\leq \frac{k}{k+1}h_K(u_{k+1})-\frac{1}{k+1}h_K(-u_{k+1}).
    \end{aligned}
    \end{equation}
    Similarly, if $b_k>0$, the inequality (\ref{7}) always has a solution. 
    
    If $b_k=0$, consider $P_{u_{k}^{\perp}}(K)$ as a $(k-1)$-dimensional convex body and by above discussion there exist $a_1,\ldots,a_{k-1}\in \mathbb{R}$ such that
    \begin{equation}\label{8}
        \begin{aligned}
        a_1&\leq \frac{k}{k+1}h_K(u_1)-\frac{1}{k+1}h_K(-u_1),
        \\
        a_2&\leq \frac{k}{k+1}h_K(u_2)-\frac{1}{k+1}h_K(-u_2),
        \\
        &\vdots
        \\
        a_{k-1}&\leq \frac{k}{k+1}h_K(u_{k-1})-\frac{1}{k+1}h_K(-u_{k-1}),
        \\
        a_1b_1+\cdots+a_{k-1}b_{k-1}&\leq \frac{k}{k+1}h_K(u_{k+1})-\frac{1}{k+1}h_K(-u_{k+1}).
    \end{aligned}
    \end{equation}
    Besides we can choose $a_k$ small enough such that $a_k\leq \frac{k}{k+1}h_K(u_k)-\frac{1}{k+1}h_K(-u_k)$. Therefore the inequality (\ref{7}) always has a solution. 
    
    If $b_k<0$, denote $\cap_{i=1}^{k+1} H_{u_i,h_K(u_i)}^-$ by $L_{k+1}$. In particular, $L_{k+1}$ is a simplex with $K\subset L_{k+1}$ and $A_{L_{k+1}}\ne \emptyset$ according to Lemma \ref{3.3}. Moreover, 
    \begin{align*}
        h_K(-u_i)\leq h_{L_2}(-u_i)\enspace \text{and}\enspace h_K(u_i)=h_{L_2}(u_i)
    \end{align*}
    for $i=1,\ldots,k+1$. Thus $A_{L_{k+1}}\subset \cap_{i=1}^{k+1} H_{u_{i},K}^-$ and $\cap_{i=1}^{k+1} H_{u_i,K}^-\ne \emptyset$. Therefore $A_K\ne \emptyset$ and Theorem \ref{3.1} is right when $n=k$. Theorem \ref{3.1} is right by induction.
\end{proof}
Now we can prove Theorem \ref{1}.
\begin{proof}[\rm{\textbf{Proof of Theorem} \ref{1}}]
According to Theorem \ref{3.1}, there exists $a\in A_K$ and 
\begin{align*}
    V(-K[1],K[n-1])&=V(-K+a[1],K[n-1])
    \\
    &=\frac{1}{n}\int_{\s}h_{-K+a}(u)dS_K(u)
    \\
    &\leq\int_{\s}h_{K-a}(u)dS_K(u)
    \\
    &=nV(K).
\end{align*}
For the equality case, $h_{-K+a}(u)=n h_{K-a}(u)$ for every $u\in \text{supp}S_K$. Since $K$ is a convex body, there are $u_1,\ldots,u_{n+1}\in \text{supp}S_K$ such that every $n$ vectors of $\{u_1,\ldots,u_{n+1}\}$ are affinely independent. Then $h_{-K+a}(u_i)=n h_{K-a}(u_i)$  means that $a$ lies in boundary of $H_{u_i,K}^-$ for every $i=\{1,\ldots,n+1\}$, which induces that $\cap_{i=1}^{n+1} H_{u_i,K}^-$ is a one point set. Denote $\cap_{i=1}^{n+1} H_{u_i,h_K(u_i)}^-$ by $L_n$ which is a simplex. Since $A_{L_n}\subset \cap_{i=1}^{n+1} H_{u_i,K}^-$, we have
$$
h_K(-u_i)=h_{L_n}(-u_i)
$$
for $i=1,\ldots,n+1$ and every vertex of $L_n$ belongs to $K$. Moreover $K\subset L_n$ and $K=L_n$. Therefore $K$ must be a simplex when the equality holds and the equality holds when $K$ is a simplex by Lemma \ref{3.3}.

\end{proof}

\section{Another way to $-K\subset nK$}
From former sections, Theorem \ref{1} is deduced by that $-K\subset nK$. We provide a completely new proof on $-K\subset nK$. Before proving $-K\subset nK$, some essential lemmas are required.

\begin{lem}\label{4.1}
    For any positive integer $m>1$ and any concave function $f:[0,1] \to [0,\infty)$, 
    \begin{equation}\label{9}
        \int_0^1\left(r-\frac{1}{m+1}\right)f^{m-1}(r)dr\geq 0
    \end{equation}
    with equality holds if and only if $f(1)=0$ and $f$ is linear.
\end{lem}
\begin{proof}
    Let $g(r)=f(r)+\frac{m+1}{m}f\left(\frac{1}{m+1}\right)r-\frac{m+1}{m}f\left(\frac{1}{m+1}\right)$. Notice that $g\left(\frac{1}{m+1}\right)=0$, $g(1)=f(1)\geq 0$ and $g$ is concave. Thus $g(r)\leq 0$ for $0\leq r\leq \frac{1}{m+1}$ and $g(r)\geq 0$ for $\frac{1}{m+1}\leq r \leq 1$ since $g$ is concave. Therefore
    \begin{align*}
    \int_0^1\left(r-\frac{1}{m+1}\right)f^{m-1}(r)dr&\geq \int_0^1\left(r-\frac{1}{m+1}\right)\left(\frac{m+1}{m}f\left(\frac{1}{m+1}\right)-\frac{m+1}{m}f\left(\frac{1}{m+1}\right)r\right)^{m-1}dr  
    \\
    &=0.
    \end{align*}
    The equality holds if and only if $g(r)=0$ for every $r\in [0,1]$, which is equivalent to that $f(1)=0$ and $f$ is linear.
\end{proof}
Back to convex bodies, we have the famous Brunn-Minkowski inequality\cite{schneider2014convex}.
\begin{thm}[the Brunn-Minkowski inequality]
    If $K, L$ are convex bodies in $\r$, then
    $$
    V(K+L)^{\frac{1}{n}}\geq V(K)^{\frac{1}{n}}+V(L)^{\frac{1}{n}}
    $$
    with equality if and only if K and L are homothetic.
\end{thm}
The following lemma as a famous corollary of the Brunn-Minkowski inequality connects convex bodies with concave functions.
\begin{lem}\label{4.3}
    If $K$ is convex body and $L$ is a $k$-dimensional convex set in $\r$, then the function
    \begin{align*}
        g(x)=V_k(K\cap(x+L))^{\frac{1}{k}},\enspace x\in \r,
    \end{align*}
    is concave on its support, where $V_k$ denotes the $k$-dimensional volume.
\end{lem}
After all these preparations, now we can prove $-K\subset nK$.
\begin{thm}\label{4.4}
    If $K$ is a convex body in $\r$ with centroid at origin, then $-K\subset nK$.
\end{thm}
\begin{proof}
    $-K\subset nK$ is equivalent to $h_K(-u)\leq nh_K(u)$ for every $u\in \s$. By definition, 
    \begin{align*}
        \int_K xdx=0.&\Longleftrightarrow \int_{-h_K(u)}^{h_K(-u)}\int_{K\cap(-ru+u^{\perp})}y-ru d\mathcal{H}^{n-1}(y)dr=0.
        \\
        &\Longrightarrow \int_{-h_K(u)}^{h_K(-u)}rV_{n-1}(K\cap(-ru+u^{\perp})) dr=0.
    \end{align*} 
   Here we denote $\int_{-h_K(u)}^{t}V_{n-1}(K\cap(-ru+u^{\perp})) dr$ by $V(t)$. Thus 
   \begin{align*}
       \int_{-h_K(u)}^{h_K(-u)}rV(r) dr=0. &\Longleftrightarrow
       rV(r)|_{-h_K(u)}^{h_K(u)}=\int_{-h_K(u)}^{h_K(-u)}V(r) dr.
       \\
       &\Longleftrightarrow h_K(-u)V(K)=\int_{-h_K(u)}^{h_K(-u)}V(r) dr.
   \end{align*}
   Now we denote $h_K(-u)+h_K(u)$ by $w(u)$. Therefore
  \begin{align*}
      h_K(-u)\leq nh_K(u).&\Longleftrightarrow h_K(-u)\leq \frac{n}{n+1}w(u).
      \\
      &\Longleftrightarrow \int_{-h_K(u)}^{h_K(-u)}V(r) dr\leq \frac{n}{n+1}w(u)V(K).
      \\
      &\Longleftrightarrow \frac{1}{n+1}w(u)V(K)\leq \int_0^{w(u)}rV_{n-1}(K\cap(-(r-h_K(u))u+u^{\perp})) dr.
  \end{align*}
  Let $S(r)=V_{n-1}(K\cap(-(r-h_K(u))u+u^{\perp}))$ and $f(r)=S^{\frac{1}{n-1}}(r/w(u))$. We have 
  \begin{align*}
      \frac{1}{n+1}w(u)V(K)\leq \int_0^{w(u)}rS(r) dr.&\Longleftrightarrow \frac{\int_0^{w(u)}rS(r) dr}{w(u)\int_0^{w(u)}S(r) dr}\geq \frac{1}{n+1}.
      \\
      &\Longleftrightarrow \frac{\int_0^1 rf^{n-1}(r) dr}{\int_0^1 f^{n-1}(r) dr}\geq \frac{1}{n+1}.
      \\
      &\Longleftrightarrow \int_0^1\left(r-\frac{1}{n+1}\right)f^{n-1}(r)dr\geq 0.
  \end{align*}
  The above inequality holds true according to Lemma \ref{4.1} and Lemma \ref{4.3}. Thus $h_K(-u)\leq nh_K(u)$ and $-K\subset nK$.
 \end{proof}
Here we can prove Theorem \ref{1} again.
 \begin{proof}
     According to Theorem \ref{4.4}, we have 
     $$
     V(-K[1],K[n-1])\leq nV(K).
     $$
     If the equality holds, $h_K(-u)=nh_K(u)$ for every $u\in\text{supp}S_K$ when $K$'s centroid is at origin. Moreover $V_{n-1}^{\frac{1}{n-1}}(K\cap(-ru+u^{\perp}))$ is linear and $V_{n-1}(K\cap(h_K(-u)u+u^{\perp}))=0$ by Lemma \ref{4.1}. Thus
     $$
     \frac{1}{n}h_K(u)V_{n-1}(K\cap(h_K(u)u+u^{\perp}))=\frac{1}{n}h_K(u)S_K(u)=\frac{1}{n+1}V(K)
     $$
     and $\text{supp}S_K$ has precisely $n+1$ elements. Therefore $K$ must be a simplex.
 \end{proof}

\end{document}